\documentclass[11pt,sumlimits,intlimits]{amsart}
\usepackage{url,cite,hyperref,amsmath,amsthm,amssymb,mathtools}
\usepackage[margin=1.25in]{geometry}
\usepackage{xcolor}
\usepackage{comment}

\newtheorem{thmx}{Theorem}
\newtheorem{corx}[thmx]{Corollary}

\newtheorem{theorem}{Theorem}[section]
\newtheorem*{theorem*}{Theorem}

\newtheorem{corollary}[theorem]{Corollary}
\newtheorem{lemma}[theorem]{Lemma}
\theoremstyle{definition}
\newtheorem{definition}[theorem]{Definition}
\newtheorem{remark}[theorem]{Remark}

\newcommand{\la}{\langle}
\newcommand{\ra}{\rangle}
\newcommand{\Z}{\mathbb{Z}}

\newcommand{\C}{\mathbb{C}}

\newcommand{\T}{\mathbb{T}}
\newcommand{\Q}{\mathbb{Q}}

\newcommand{\ZZ}{\mathbb{Z}[\frac{1}{p}]}
\newcommand{\id}{\operatorname{id}}
\newcommand{\tr}{\operatorname{Tr}}
\newcommand{\GL}{\operatorname{GL}}
\newcommand{\e}{\mathbf{e}}

\numberwithin{equation}{section}

\title{Residually finite amenable groups that are not Hilbert-Schmidt Stable}
\author{Caleb Eckhardt}
\address{Department of Mathematics, Miami University, Oxford, OH, 45056}
\email{eckharc@miamioh.edu}

\begin{document}
\thanks{This work was partially supported by an AMS-Simons research enhancement grant for PUI faculty.}
\maketitle
\begin{abstract} We construct the first examples of residually finite amenable groups that are not Hilbert-Schmidt (HS) stable. We construct finitely generated, class 3 nilpotent by cyclic examples and solvable linear finitely presented examples.  This also provides the first examples of amenable groups that are very flexibly HS-stable but not flexibly HS-stable and the first examples of residually finite amenable groups that are not locally HS-stable. Along the way we exhibit (necessarily not-finitely-generated) class 2 nilpotent groups  $G = A\rtimes \Z$ with $A$ abelian such that the periodic points of the dual action are dense but it does not admit dense periodic measures.  Finally we use the Tikuisis-White-Winter theorem to show all of the examples are not even operator-HS-stable; they admit operator norm almost homomorphisms that can not be HS-perturbed to true homomorphisms.
\end{abstract}
\section{Introduction}   
This paper fills in a few gaps in the literature relating (local) Hilbert-Schmidt stability to residual finiteness of amenable groups. There has been renewed interest recently in stability questions due in part to their connection with strategies for locating non-sofic or non-hyperlinear groups.  We refer the reader to \cite{Becker19,Becker20} for more on this motivation. In particular it is desirable to establish criteria for a group to be HS-stable (or not).  In this direction we show that residual finiteness does not guarantee HS-stability for finitely presented amenable groups.

 Let $U(n)$ denote the group of unitary $n\times n$ matrices equipped with the metric $d_n(u,v) = \tr_n((u-v)^*(u-v))^\frac{1}{2}$ where $\tr_n$ is the normalized trace (i.e. $\tr_n$ of the identity matrix is 1). A group $G$ is \textbf{Hilbert-Schmidt (HS) stable} if for every sequence of functions $\pi_k:G \to U(n_k)$ satisfying $\lim_kd_{n_k}(\pi_k(gh),\pi_k(g)\pi_k(h))= 0$ for all $g,h\in G$ there is a sequence of homomorphisms $\sigma_k:G\to U(n_k)$ such that $\lim_kd_{n_k}(\sigma_k(g),\pi_k(g)) = 0$ for all $g\in G.$

The F\o lner condition shows that every finitely generated HS-stable amenable group is residually finite. We show the converse does not hold in general answering \cite[Question 1.1]{Eckhardt23}
\begin{thmx} \label{thm:A}  There are infinitely many non-isomorphic finitely generated, residually finite, non HS-stable groups of the form $N\rtimes \Z$ where $N$ is class 3 nilpotent.  Moreover these groups are not operator-HS-stable (See Theorem \ref{thm:sform}).
\end{thmx}
All finitely generated nilpotent-by-finite groups are HS-stable by \cite{Levit24}. In this sense the nilpotent-by-cyclic examples of Theorem \ref{thm:A} are the simplest one could hope for.  
Becker and Lubotzky  \cite{Becker20} and Ioana \cite{Ioana20a} gave examples of non-amenable residually finite groups that are not HS-stable, our contribution provides amenable examples. We also mention that Hadwin and Shulman in \cite[Example 3.11]{Hadwin18a} 
showed a C*-algebraic version of Theorem \ref{thm:A} by
constructing an amenable C*-algebra that is residually finite-dimensional but is not matricially stable.

Local Hilbert-Schmidt stability was recently introduced and systematically studied by Fournier-Facio, Gerasimova and Spaas \cite{Fournier24}. Every HS-stable group is trivially locally HS-stable and the notions coincide for finitely presented groups. We modify the examples of Theorem \ref{thm:A} to strengthen its conclusion

\begin{thmx} \label{thm:fpres} There are infinitely many mutually non-isomorphic  finitely presented solvable subgroups of $GL(5,\mathbb{Q})$ that are residually finite and not locally HS-stable (and hence not HS-stable).
\end{thmx}

Recent works on HS-stability \cite{Eckhardt23,Levit24} realize the importance of the relationship between dense periodic points and dense periodic measures. We describe a simple case and refer the reader to \cite{Levit24} for more information.  Let $A$ be an abelian group and $\alpha$ an automorphism of $A.$  We say that the dual dynamical system  $(\widehat{A},\widehat\alpha)$ admits \textbf{dense periodic measures} if the set of all $\widehat\alpha$-invariant probability measures on $\widehat A$ with finite support is weak*-dense in the set of all $\widehat\alpha$-invariant probability measures on $\widehat A$.  By \cite[Proposition 4.8, Corollary 5.19]{Eckhardt23} the group $A\rtimes \Z$ is a union of residually finite groups precisely when the periodic points of $(\widehat{A},\widehat\alpha)$ are dense and the group is HS-stable precisely when $(\widehat{A},\widehat\alpha)$ admits dense periodic measures.

While proving Theorem \ref{thm:A} we show the following answering \cite[Question 1.2]{Eckhardt23}
\begin{thmx} \label{thm:C} Let $p$ be prime and $G_p = \ZZ^2\rtimes \Z$ where the automorphism is given by $\alpha = \left[  \begin{array}{cc} 1 & 1 \\ 0 & 1 \end{array} \right].$ Then $G_p$ is residually finite but not HS-stable.
Hence the periodic points of $(\widehat{\ZZ}^2,\widehat\alpha)$ are dense but the system does not admit dense periodic measures.
\end{thmx}
We remark that the groups $G_p$ are not finitely generated. It appears open \cite[Question 2]{Levit24} whether or not finitely generated groups of the form $A \rtimes \Z$ with $A$ abelian always admit dense periodic measures.

(Very) flexible notions of HS-stability were introduced by Becker and Lubotzky \cite{Becker20}.
Roughly, a group is very flexibly HS-stable if the homomorphisms $\sigma_n$ are allowed to take values in a larger unitary group and a cut down of $\sigma_n$ approximates $\pi_n.$ 

It was shown in \cite[Section 6]{Eckhardt23} that a finitely generated amenable group is flexibly HS-stable if and only if it is HS-stable \footnote{It is apparently open whether or not these notions coincide in general} and very flexibly HS-stable if and only if it is residually finite.  Therefore Theorem \ref{thm:A} immediately yields

\begin{corx} \label{thm:B}  The groups of Theorems \ref{thm:A},\ref{thm:fpres} and \ref{thm:C} are very flexibly HS-stable but not flexibly HS-stable.
\end{corx}

Finally, we describe our strategy for constructing the examples.  For \emph{amenable} groups, Hadwin and Shulman rephrased HS-stability in terms of tracial approximations.  A \textbf{trace} $\tau:G\to \C$ is a normalized positive definite function that is constant on conjugacy classes. We call $\tau$ \textbf{finite-dimensional} if there is a finite-dimensional unitary representation $\pi:G\to U(n)$ such that $\tau = \tr_n\circ \pi.$ Then \cite[Theorem 4]{Hadwin18} shows an amenable group $G$ is HS-stable if and only if every trace on $G$ is the pointwise limit of finite-dimensional traces.

The idea is to take a residually finite group $G$ and a normal subgroup $N$ such that $G/N$ is not HS-stable and then show that the trace on $G$ induced from the canonical trace on $G/N$ is not the pointwise limit of finite-dimensional traces.  

Of course some care must be taken in the choice of groups $G,N.$  Indeed it was shown in \cite[Corollary 5.12]{Eckhardt23} that there is an HS-stable group $G$ and a central subgroup $N\leq G$ such that $G/N$ is not HS-stable. The centrality of $N$ was a crucial fact in showing HS-stability
 in that case.  So we start with pairs $G,N$ with $G$ residually finite and $G/N$ not residually finite so the action of $G$ on $N$ is non-trivial, but trivial enough so we can keep track of all of the finite-dimensional traces.  
 
This strategy leads to the groups of Theorem \ref{thm:C}. Notice that $\ZZ$ is residually finite but any homomorphism $\pi$ of $\ZZ$ into a finite group with $1\in\text{ker}(\pi)$ is trivial (in particular the quotient $\ZZ/\Z$ is not residually finite). We introduce an action on $\ZZ^2$ that leaves $\Z^2$ invariant (thus producing a trace) while severely restricting the finite-dimensional representations of $\ZZ^2\rtimes \Z$ (and hence the finite-dimensional traces) in a similar way that the finite group quotients of $\ZZ$ are restricted. This restriction shows the trace on $\ZZ^2\rtimes \Z$ induced from the canonical trace on $(\ZZ/\Z)^2\rtimes \Z$ is not the pointwise limit of finite-dimensional traces. 

Since these groups are not finitely generated we embed them into a finitely generated group without destroying the obstruction to HS-stability thus proving Theorem \ref{thm:A}. The groups of Theorem \ref{thm:A} are not finitely presented so we adapt the techniques of Abels \cite{Abels79} to produce the finitely presented examples of Theorem \ref{thm:fpres}.
 
We finish the paper by showing how the Tikuisis-White-Winter theorem and Theorems \ref{thm:A},\ref{thm:fpres} and \ref{thm:C} combine to show these groups enjoy a strong negation of HS-stability (Theorem \ref{thm:sform}).

\subsection*{Acknowledgements} I thank Francesco Fournier-Facio for making me aware of \cite{Fournier24} and local HS-stability in general, which provided extra motivation to find finitely presented examples.  I also thank Dan Farley for a helpful conversation about groups.
 
\section{Examples}
For an abelian group $A$ we write $\widehat{A}$ for the Pontraygin dual group.
Let $p$ be prime and let $\ZZ = \la  p^n:n<0 \ra \leq (\Q,+).$ 
This lemma shows the simple algebraic fact we exploit to produce all of the examples.
\begin{lemma} \label{lem:finorder} Let $m\geq 1$ and suppose that $\chi \in \widehat {\ZZ}$ has finite order $m.$ Then $p$ does not divide $m$. In particular if $\chi(1)=1$, then $\chi$ is the trivial character.
\end{lemma}
\begin{proof} By assumption there is some $k\in \Z$ and $n\geq 0$ such that $\chi(\frac{k}{p^n})$ has order $m.$
If $p$ divides $m$, then
\begin{equation*}
1 = \chi\left(\frac{k}{p^{n+1}}\right)^m = \left(  \chi\left(\frac{k}{p^{n+1}}\right)^p \right)^{\frac{m}{p}} = \chi\left(\frac{k}{p^n}\right)^{\frac{m}{p}} 
\end{equation*}
contradicting the fact that $\chi(\frac{k}{p^n})$ has order $p.$

For the last claim, suppose that $\chi(1)=1.$  Then for each integer $n\geq 1$, we have that $\chi(p^{-n})$ is an $m$th root of unity and a $p^n$th root of unity. Since $p$ does not divide $m$ we have $\chi(p^{-n})=1.$
\end{proof}

\subsection{Residually finite, infinitely generated non-HS stable groups}

\begin{definition}\label{def:G} Consider the automorphism $\alpha = \left[  \begin{array}{cc} 1 & 1 \\ 0 & 1  \end{array} \right]$ of $\ZZ^2.$  Define 
\newline
$G_p = \ZZ^2 \rtimes_\alpha \Z.$ We follow standard convention and write elements in $G_p$ as triples $(a,b,c)$ where $a,b\in \ZZ$, $c\in \Z$ and 
\begin{equation*}
(a_1,b_1,c_1)(a_2,b_2,c_2) = (a_1+a_2+c_1b_2,b_1+b_2,c_1+c_2).
\end{equation*}
\end{definition}
\begin{lemma} \label{lem:nilrf} The group $G_p$ is class 2 nilpotent and residually finite.
\end{lemma}
\begin{proof} The Heisenberg group $\mathbb{H}_3(\ZZ)$ is class 2 nilpotent and residually finite. Note that the map $(a,b,c) \to \left[  \begin{array}{ccc} 1 & b & a-bc \\ 0 & 1 & -c\\ 0 & 0 & 1   \end{array} \right]$ is an injective group homomorphism from $G_p$ into $\mathbb{H}_3(\ZZ).$ 
\end{proof}
\begin{remark} \label{rem:Heis} The Heisenberg group $\mathbb{H}_3(\ZZ)$ is HS-stable by \cite[Theorem 5.11]{Eckhardt23} and we will show that its subgroup $G_p$ is not HS-stable.  The difference is that the image of $\Z^2$ in $\mathbb{H}_3(\ZZ)$ is not a normal subgroup, hence the ``bad" trace that $G_p$ admits does not extend to a trace on $\mathbb{H}_3(\ZZ).$
\end{remark}
We now show that $G_p$ is not HS-stable.  The finite dimensional irreducible representations of $G_p$ are more-or-less characterized by the periodic points of the dual dynamical system $(\widehat{\ZZ}^2,\widehat \alpha)$ so we first characterize those.  Recall that for every character $\chi$ of $\ZZ^2$ there are $\chi_1,\chi_2\in \widehat{\ZZ}$ such that $\chi(a,b) = \chi_1(a)\chi_2(b).$

\begin{lemma} \label{lem:ppoints} Every $m$-periodic point of the dynamical system $(\widehat{\ZZ}^2,\widehat \alpha)$ is of the form $\chi = (\chi_1,\chi_2)$ where $\chi_2$ is any element of $\widehat{\ZZ}$ and $\chi_1\in \widehat{\ZZ}$ has order $m.$ Moreover $p$ does not divide $m$.
\end{lemma}
\begin{proof} Suppose $\chi=(\chi_1,\chi_2)$ has period $m$.  Then for all $(a,b)\in \ZZ^2$ we have
\begin{equation*}
\chi_1(a)\chi_2(b) = \chi(a,b) = \chi(\alpha^m(a,b)) = \chi(a+mb,b) = \chi_1(a)\chi_2(b)\chi_1(b)^m.
\end{equation*}
Hence $\chi_1$ has order $m.$ The fact that $p$ does not divide $m$ follows from Lemma  \ref{lem:finorder}.
\end{proof}
\begin{definition} We let $\tr_n$ denote the normalized trace on $n\times n$ complex matrices.
\end{definition}
\begin{lemma} \label{lem:fdtraces}  Let $\pi$ be an irreducible $n$-dimensional representation of $G_p.$ Then $p$ does not divide $n$ and there are characters $\chi_1,\chi_2 \in \widehat{\ZZ}$ and $\phi\in\widehat{\Z}$ such that $\chi_1$ has period $n$ and
\begin{equation*}
\tr_n\circ\pi(a,b,c) = \left\{ \begin{array}{ll} \chi_1(a)\chi_2(b)\phi(c) & \text{ if }   \chi_1(b) = 1 \text{ and } n|c\\
                                                                       0 & \text{ otherwise } \end{array} \right.
\end{equation*}
\end{lemma}
\begin{proof} Since $\pi$ is irreducible there is a $\widehat{\alpha}$-periodic point $\chi = (\chi_1,\chi_2)$ with period $n$ such that (up to unitary equivalence) we have
\begin{equation*}
\pi(a,b,0) = \bigoplus_{i=0}^{n-1} \chi(\alpha^i(a,b)) = \bigoplus_{i=0}^{n-1} \chi(a+bi,b) = \bigoplus_{i=0}^{n-1} \chi_1(a)\chi_2(b)\chi_1(b)^i
\end{equation*}
and there is some $\lambda\in\T$ such that
\begin{equation*}
\pi(0,0,1) = \left[ \begin{array}{ccccc} 0 & \lambda & 0 & \cdots & 0\\
                                                         0 & 0 & \lambda & \cdots & 0\\
                                                         \vdots & \vdots & \vdots & \ddots & \vdots\\
                                                         0 & 0 & 0 & \cdots & \lambda\\
                                                         
                                                         \lambda & 0 & 0 & \cdots & 0    \end{array}   \right].
\end{equation*}
Define $\phi \in \widehat \Z$ by $\phi(1)=\lambda.$ 

Notice that if $n$ does not divide $c$ then $\pi(a,b,c)$ has all diagonal entries equal to 0 hence $\tr_n\circ \pi(a,b,c)=0.$ 
 On the other hand if $n|c$ then $\pi(0,0,c) = \phi(c)\id,$ hence
\begin{equation*}
\tr_n\circ\pi(a,b,c) = \tr_n\left(\phi(c)\bigoplus_{i=0}^{n-1} \chi_1(a)\chi_2(b)\chi_1(b)^i\right) = \frac{1}{n}\chi_1(a)\chi_2(b)\phi(c)\sum_{i=0}^{n-1} \chi_1(b)^i.
\end{equation*}
If $\chi_1(b)=1$ we have $\tr_n\circ\pi(a,b,c)=\chi_1(a)\chi_2(b)\phi(c).$   By Lemma \ref{lem:ppoints} $\chi_1$ has order $n.$ Hence if $\chi_1(b)\neq 1$, then $\chi_1(b)$ is a non-trivial $n$th root of unity, thus $\tr_n\circ\pi(a,b,c)=0.$
\end{proof}
Theorem \ref{thm:C} is an immediate consequence of the following theorem and \cite[Theorem 6.2, Corollary 6.6]{Eckhardt23}
\begin{theorem} \label{thm:main} Let $H_p = \la  (1,0,0),(0,1,0) \ra \leq G_p.$
Suppose that $\tau_n:G_p \to \C$ is a sequence of finite-dimensional traces such that $\lim_n \tau_n(0,1,0)= 1.$  Then $\lim_n \tau_n(a,0,0) =  1$ for every $a\in \ZZ.$ In particular the trace defined by
\begin{equation*}
\tau(g) = \left\{ \begin{array}{cc} 1 & \textup{ if }g\in H_p\\
                                             0 & \textup{ if }g\not\in H_p   \end{array} \right.
 \end{equation*}
is not the pointwise limit of finite-dimensional traces.  Hence $G_p$ is not HS-stable.  Moreover if $G_p$ embeds into an amenable group $K$ such that the image of $H_p$ is normal in $K$, then $K$ is also not HS-stable.
\end{theorem}
\begin{proof} Let $a\in\ZZ.$ By Lemma \ref{lem:fdtraces} the extreme finite-dimensional traces of $G_p$ are parametrized by the set
\begin{equation*}
X = \{ (\chi_1,\chi_2,\phi) \in \widehat{\Z[\tfrac{1}{p}]}\times \widehat{\Z[\tfrac{1}{p}]} \times \widehat{\Z} : \chi_1 \text{ has finite order }  \}.
\end{equation*}
For each $x=(\chi_1,\chi_2,\phi)\in X$, let $\tau_x$ be the corresponding trace from Lemma \ref{lem:fdtraces}. We say that $\tau_x$ is of type I if $\chi_1(1)=1$ and of type 2 if $\chi_1(1)\neq 1.$  By Lemma \ref{lem:finorder} if $\tau_x$ is a  type 1 trace, then $\chi_1$ is the trivial character. In particular for every type 1 trace $\tau_x$ by Lemma \ref{lem:fdtraces} we have
\begin{equation}\label{eq:t1tracetriv}
\tau_x(a,0,0) = \chi_1(a) = 1 \text{ for all type 1 traces }\tau_x.
\end{equation}

Every finite-dimensional trace is a convex combination of extreme finite-dimensional traces.  Therefore for each $n$ there are scalars $0\leq \mu_n\leq 1$ and traces $\tau_{1,n}$ and $\tau_{2,n}$ such that $\tau_{1,n}$ is a convex combination of type 1 traces, and
$\tau_{2,n}$ is a convex combination of type 2 traces and $\tau_n = \mu_n\tau_{1,n} + (1-\mu_n)\tau_{2,n}.$

By Lemma \ref{lem:fdtraces} we have $\tau_{2,n}(0,1,0)=0$, so
\begin{equation*}
\tau_n(0,1,0) = \mu_n\tau_{1,n}(0,1,0) + (1-\mu_n)\tau_{2,n}(0,1,0) = \mu_n\tau_{1,n}(0,1,0).
\end{equation*}
Since $\tau_n(0,1,0)$ converges to 1 we have $\mu_n$ also converges to 1. 
\newline
Since $\tau_{1,n}$ is a linear combination of type 1 traces, by (\ref{eq:t1tracetriv}) we have $\tau_{1,n}(a,0,0)=1$, hence
\begin{equation*}
\tau_n(a,0,0) =  \mu_n + (1-\mu_n)\tau_{2,n}(a,0,0).
\end{equation*}
Since $\mu_n$ converges to 1 we have that $\tau_n(a,0,0)$ also converges to 1. That $G_p$ is not HS-stable now follows from \cite[Theorem 4]{Hadwin18} (see also Section \ref{rem:apphomo} for a very direct proof adapted to our situation).

For the final claim suppose that $H_p\leq G_p \leq K$ with $H_p$ normal in $K.$  Since $H_p$ is normal the trace $\tau$ on $G_p$ extends to $K$ by setting $\tau(g) = 0$ for all $g\not\in G_p.$ If $K$ is HS-stable then there are finite-dimensional representations $\sigma_k:K \to U(n_k)$ such that $\tau$ is the pointwise limit of $\tr_{n_k} \circ \sigma_k.$ Then $\tau|_{G_p}$ is the pointwise limit of $\tr_{n_k} \circ (\sigma_k|_{G_p}),$ a contradiction. 
\end{proof}
Combining Lemma \ref{lem:nilrf}, Remark \ref{rem:Heis} and the previous theorem we make the following observation.
\begin{corollary} The non-HS stable group $G_p$ embeds as a normal subgroup inside of the HS-stable Heisenberg group $\mathbb{H}_3(\ZZ).$  In particular, HS-stability for amenable groups is not necessarily preserved by taking normal subgroups.
\end{corollary}

\subsection{Residually finite, finitely generated non HS-stable groups}

Next we use standard ideas to construct a finitely generated example that contains $G_p$ and keeps the HS-stability obstruction. These groups were the simplest we could think of but they are not finitely presented by \cite[Corollary C]{Bieri78} (it is easy to see $K_p$ violates condition iv.c).  In the next section we construct more complicated finitely presented examples.

\begin{definition} \label{def:K}  Let $K_p$ be the subgroup of upper triangular elements of $\GL(5,\ZZ)$ defined by
\begin{equation*}
K_p = \left\{   \left[  \begin{array}{ccccc}  1 & x_{12} & 0 & x_{14} & x_{15}\\
                                                          & p^n & 0 & x_{24} & x_{25} \\ 
                                                          && 1 & 0 & 0\\
                                                          &&& 1 & m\\
                                                          &&&& 1\end{array}  \right]  : x_{ij} \in \Z[\tfrac{1}{p}] \text{ and }n,m\in\Z   \right\}.
\end{equation*}
Let $N\leq K_p$ be the subgroup consisting of all elements with $n=0.$  Let $ G_K \leq N$ be the subgroup that further has $x_{12}=x_{24}=x_{25}=0.$ Let $H_K\leq G_K$ be the subgroup that further has $m=0$ and $x_{14},x_{15}\in \Z.$
\end{definition}
\begin{proof}[Proof of Theorem \ref{thm:A}]  Since $\GL(5,\ZZ)$ is residually finite so is $K_p.$ The subgroup $N$ is readily seen to be class 3 nilpotent. Let $a = \text{diag}(1,p,1,1,1)\in K_p.$  Then one easily sees that $K_p\cong N\rtimes \Z$ where the action is conjugation by $a.$

Next we show $K_p$ is finitely generated.
For $1\leq i,j \leq 5$ let  $\e_{ij}$ be the $5\times 5$ matrix with the $i,j$ entry equal to 1 and all other entries equal to 0.  Let $1\in \GL(5,\ZZ)$ denote the identity. We claim that
\begin{equation} \label{eq:Kgens}
K_p = \la  a, 1+\e_{12},1+\e_{24},1+\e_{25},1+\e_{45}  \ra.
\end{equation}
Let $g = \left[  \begin{array}{ccccc}  1 & x_{12} & 0 & x_{14} & x_{15}\\
                                                          & 1 & 0 & x_{24} & x_{25} \\ 
                                                          && 1 & 0 & 0\\
                                                          &&& 1 & m\\
                                                          &&&& 1\end{array}  \right] \in N.$ Notice that
                                                          \begin{equation}\label{eq:Nelts}
g = (1+m\e_{45})(1+x_{24}\e_{24})(1+x_{25}\e_{25})(1+x_{12}\e_{12})(1+x_{14}\e_{14})(1+x_{15}\e_{15}).
\end{equation}

 Suppose that $i<j$ and $r<s$ and $x,y\in \ZZ.$  Then (the brackets indicate \emph{group} commutators: $[g,h] = g^{-1}h^{-1}gh$):
 \begin{equation}\label{eq:commrel}
 (1+x\e_{ij})^{-1} = 1-x\e_{ij}\quad \text{ and } \quad 
                             [1+x\e_{ij},1+y\e_{rs}] = \left\{ \begin{array}{cc} 1+xy\e_{is} & \text{ if } j=r\\
                                                                                                        1-xy\e_{rj} & \text{ if }i=s\\
                                                                        1 & \text{ else }  \end{array} \right.
                                                                         \end{equation}

For the elements $1+\e_{ij} \in \{1+\e_{12},1+\e_{24},1+\e_{25}\}$ we have
\begin{equation}\label{eq:allent}
\la a^n(1+\e_{ij})a^{-n} : n\in \Z  \ra = \{ 1+x\e_{ij} : x\in \Z[\tfrac{1}{p}]  \}.
\end{equation}
We emphasize that (\ref{eq:allent})  does not hold for $1+\e_{45}.$  

By (\ref{eq:commrel}) we have
\begin{equation*}
[1+x\e_{12},1+\e_{24}] = 1 + x\e_{14} \quad \text{for all }x\in \Z[\tfrac{1}{p}]
\end{equation*}
and 
\begin{equation*}
[1+x\e_{14},1+\e_{45}] = 1 + x\e_{15} \quad \text{for all }x\in \Z[\tfrac{1}{p}].
\end{equation*}
Therefore by (\ref{eq:Nelts}) we obtain (\ref{eq:Kgens}). This shows $K_p$ is finitely generated.

By (\ref{eq:commrel}) the elements $1+ \e_{12},1+\e_{24},1+\e_{25}$ all commute with every element of $H_K.$ One also checks that $a$ commutes with $H_K$. Moreover by (\ref{eq:commrel}) we have that $1+\e_{45}$ commutes with $1+\e_{15}$ and
\begin{equation*}
(1+\e_{45})(1-\e_{14})(1+\e_{45})^{-1} =  (1-\e_{14})(1+\e_{15}).
\end{equation*}
Therefore $H_K$ is a normal subgroup of $K_p.$
Moreover the conjugation above shows that $G_K$ is isomorphic to $G_p$.  Therefore $K_p$ is not HS-stable by Theorem 
\ref{thm:main}.  To see that the groups $K_p$ are pairwise non-isomorphic note that the center of $K_p$ is isomorphic to $\ZZ.$
\end{proof}

\subsection{Residually finite, finitely presented non (locally) HS-stable groups}

Now we use the ideas of \cite{Abels79} to embed the examples of the previous section into a finitely presented solvable linear group and keep the HS-stability obstruction.  To those readers familiar with the ideas in \cite{Abels79} the proof is probably self-evident.  Nonetheless we provide all of the details for the convenience of the other readers.

\begin{lemma}\label{lem:fplem}
For each prime $p$ define
\begin{equation*} 
\tilde{G}_p = \left\{ \left[  \begin{array}{ccccc}  1 & x_{12} & x_{13} & x_{14} & x_{15}\\
                                            0 & p^n & x_{23} & x_{24} & x_{25} \\ 
                                            0 & 0 & p^k & x_{34} & x_{35} \\
                                            0 & 0 & 0 & 1 & m\\
                                            0 & 0 & 0 & 0 & 1\\    \end{array}  \right]   : x_{ij} \in \Z[\tfrac{1}{p}] \text{ and }k,n,m\in\Z \right\}
\end{equation*}
Then $\tilde{G}_p$ is finitely presented.
\end{lemma}
\begin{proof}  Let $p$ be a prime and $G\leq \tilde G_p$ be the group which further satisfies $m=0.$ Since $\tilde{G}_p \cong G\rtimes \Z$, the proof will be complete once we show $G$ is finitely presented. 

Let $H$ be the group with generators $d_2,d_3, e_{ij},$ for $1\leq i\leq3$ and $i<j\leq 5$ subject to the relations
\begin{flalign*}
d_2\textbf{ relations: }   [d_2,d_3] = 1; \hspace{0.2in} d_2e_{2j}d_2^{-1} = e_{2j}^p \text{ for }j=3,4,5; \hspace{0.2in} d_2^{-1}e_{12}d_2 = e_{12}^p; &&
\end{flalign*}  
\begin{flalign*}
[d_2,e_{ij}] = 1 \text{ for all pairs }(i,j)\text{ not listed in the previous line.}
\end{flalign*}
\begin{flalign*}
d_3 \textbf{ relations: } d_3x_{3j}d_3^{-1} = e_{3j}^p \text{ for }j=4,5;\hspace{0.2in} d_3^{-1}e_{i3}d_3 = e_{i3}^p \text{ for }i=1,2 &&
\end{flalign*}
\begin{equation*}
[d_3,e_{ij}] = 1 \text{ for all pairs }(i,j)\text{ not listed in the previous line.}
\end{equation*}
\begin{flalign*}
e_{ij}\textbf{ relations: }\text{For } j\leq k, \hspace{0.2in}   [e_{ij},e_{k\ell}] = \left\{ \begin{array}{ll} 1 & \text{if }j\neq k\\
e_{i\ell} &\text{if } j=k \end{array} \right. &&
\end{flalign*}
We show that $G$ is isomorphic to $H.$  

By the main result of \cite{Abels79} (we have duplicated his notation) the relations above that involve $d_2,d_3$ and $e_{ij}$ with $j\neq 5$ are a presentation for the group $G_1\leq G$ that have $x_{i5}=0$ for $i=1,2,3.$  In fact the map  $d_2 \mapsto \text{diag}(1,p,1,1,1), d_3\mapsto \text{diag}(1,1,p,1,1)$ and $e_{ij}\mapsto 1 +\e_{ij}$ extends to an explicit isomorphism. We freely use any relation that holds in $G_1$ inside of $H.$

Moreover as in the proof of the main result of \cite{Abels79} one applies his Lemma twice to see that the relations involving $d_2,d_3$ and $e_{ij}$ with $i\neq 1$ are a presentation for the group $G_2\leq G$ that have $x_{1j}=0$ for $j=2,3,4,5.$  As above the map  $d_2 \mapsto \text{diag}(1,p,1,1,1), d_3\mapsto \text{diag}(1,1,p,1,1)$ and $e_{ij}\mapsto 1 +\e_{ij}$ extends to an isomorphism. We freely use any relation that holds in $G_2$ inside of $H.$

Let $r = \frac{k}{p^n} \in \ZZ$ with $p \not| k.$  Define $re_{12} = d_2e_{12}^nd_2^{-1}e_{12}^k.$  We similarly define $re_{ij}$ for all pairs $(i,j) \neq (1,4),(1,5).$ Note that $re_{23}$ is unambiguous since $d_2^{-1}e_{23}d_2= d_3e_{23}d_3^{-1}.$ We define $re_{14} = [re_{12},e_{24}].$ 
By the above remarks about $G_1$ and $G_2$ we have for $r,s \in \ZZ$ and $j\leq k$
\begin{equation} \label{eq:updatedcomm}
  [re_{ij},se_{k\ell}] = \left\{ \begin{array}{ll} 1 & \text{if }j\neq k\\
rse_{i\ell} &\text{if } j=k \end{array} \right. 
\end{equation}
whenever either $i\neq 1 \neq k$ or $j\neq 5 \neq \ell.$  Next we extend the above equations to the missing pairs.  Let $r,s\in \ZZ.$  Notice that $[re_{35},se_{12}]=1$ since $[d_2,d_3] = [e_{35},e_{12}]=[d_2,e_{35}]=[d_3,e_{12}]=1.$ It then follows that $[re_{35},se_{14}]=1$ since $se_{14} \in \la  se_{12},e_{24} \ra.$  Similar arguments show $[re_{25},se_{13}] = [re_{25},se_{14}] = 1.$

We need to unambiguously define $re_{15}$ and show that it is central.  This follows from the tiniest of modifications to the proof of \cite{Abels79}. 
For any group elements $x,y$ let $x^y = y^{-1}xy.$  We then have the following relationship for any elements $x,y,z$ in a group:
\begin{equation}\label{eq:Hall}
[x^y,[y,z]]\cdot[y^z,[z,x]]\cdot[z^x,[x,y]]=1.
\end{equation}
Applying this equation to $x=re_{12}, y = e_{23}$ and $z= ye_{35}$ and using the above established commutation identities we obtain  $[re_{12},se_{25}] = [re_{13},se_{35}]$ which we use to define $re_{15} = [re_{12},e_{25}].$ Finally we show $re_{15}$ is central.  Since $re_{15} \in \la re_{12},e_{25}  \ra \cap \la re_{13},e_{35}  \ra$, we have $[re_{15},se_{ij}] =1 $ for all $(i,j) \neq (2,3)$ and $[re_{15},d_j] = 1$ for $j=2,3.$

Then applying (\ref{eq:Hall}) to $x = re_{13}, y=e_{35}$ and $z=e_{23}$ we obtain $[re_{15},e_{23}]=1.$  Since $[re_{15},d_2] = 1$ we also have $[re_{15},se_{23}]=1$ showing that $re_{15}$ is central.

It now follows that we have a surjective homomorphism $\phi:H \to G$ defined by $\phi(re_{ij}) = 1 +r\e_{ij}$ and $\phi(d_i)$ sent to the corresponding diagonal elements. The above computations show that each element of $H$ can be written as
\begin{equation*}
(r_{15}e_{15})(r_{25}e_{25})(r_{35}e_{35})(r_{14}e_{14})(r_{24}e_{24})(r_{34}e_{34})(r_{13}e_{13})(r_{23}e_{23})(r_{12}e_{12})d_2^nd_3^m
\end{equation*}
for some $r_{ij}\in \ZZ$ and $n,m\in \Z.$  Since every element in $G$ can be uniquely written as above (with the obvious label changes) it follows that $\phi$ is an isomorphism.
\end{proof}

\begin{proof}[Proof of Theorem \ref{thm:fpres}] As in the proof of Theorem \ref{thm:A} we apply Theorem \ref{thm:main} to see that $\tilde{G}_p$ is not HS-stable.  By Theorem \ref{lem:fplem} the group $\tilde{G}_p$ is finitely presented, hence by \cite[Lemma 3.11]{Fournier24} it follows that $\tilde{G}$ is not locally HS-stable.
\end{proof}

\subsection{Explicit non-perturbable approximate homomorphisms} \label{rem:apphomo}

Suppose that $\Gamma$ is amenable and $\Delta \trianglelefteq \Gamma$ is a normal subgroup such that the trace $\tau$ defined by $\tau(g) = 1$ if $g\in \Delta$ and $\tau(g)=0$ if $g\not\in \Delta$ is not the limit of finite dimensional traces. For example, any of the examples of the previous sections.  Then \cite[Theorem 4]{Hadwin18} shows $\Gamma$ is not HS-stable.  Due to the mundane nature of $\tau$ it is easy to construct approximate homomorphims of $\Gamma$ that are not perturbable thus giving a direct proof of non-stability.

Let $F\subseteq \Gamma/\Delta$ be a finite subset.  For each $x\in \Gamma $ let $x_F = P_F\lambda_x P_F$ where $P_F$ is the projection of $\ell^2 \Gamma/\Delta$ onto $\ell^2F$ and $\lambda_x$ is the left regular representation.  Let $\pi_F(x)$ be any unitary in $B(\ell^2F)$ satisfying $x_F = U|x_F|.$  Then $\| \pi_F(x) - x_F  \|_2 \leq \| 1-|x_F| \|_2.$ 

Now let $(F_n)$ be a F\o lner system for $\Gamma/\Delta.$ The above observations show that $x \mapsto \pi_{F_n}(x\Delta)$ defines a HS-approximate homomorphism from $\Gamma$ into the unitary matrices with HS metric. One easily checks (as in the proof of Theorem \ref{thm:sform} below) that if this approximate homomorphism could be perturbed to an actual homomorphism then $\tau$ would be the limit of finite-dimensional traces.

\subsection{A strong negation of HS-stability} 
Finally we point out how the Tikuisis-White-Winter theorem combines with Theorem \ref{thm:main} to build approximate homomrphisms into the unitary matrices with the operator norm that can not be Hilbert-Schmidt perturbed to a sequence of true homomorphisms.  Since the operator norm dominates the Hilbert-Schmidt norm this is formally stronger than the negation of Hilbert-Schmidt stability.
 \begin{definition} Let $M_n(\C)$ be the $n\times n$ complex matrices.  Let $A\in M_n(\C)$. We let $\| A \|_\infty$ denote the operator norm and  $\| A \|_2 = \tr_n(A^*A)^\frac{1}{2}$. 
 \end{definition}
 We isolate an easy application of the spectral theorem. 
 \begin{lemma} \label{lem:basicalc} Let $A\in M_n(\C)$ be invertible with $\| A \|_\infty \leq 1.$  Let $U$ be the unitary matrix in the polar decomposition of $A.$  Then $\| U-A \|_\infty \leq \|  1-A^*A \|_\infty.$
 \end{lemma}
 
\begin{theorem} \label{thm:sform} 
Suppose that $\Gamma$ is amenable and $\Delta \trianglelefteq \Gamma$ is a normal subgroup such that the trace $\tau$ defined by $\tau(g) = 1$ if $g\in \Delta$ and $\tau(g)=0$ if $g\not\in \Delta$ is not the limit of finite dimensional traces.  There is a sequence of maps $\pi_k:\Gamma\to U(n_k) $ such that 
\begin{equation} \label{eq:almostmap}
\lim_{k\rightarrow\infty} \|  \pi_k(gh)-\pi_k(g)\pi_k(h) \|_\infty =0 \quad \text{ for all }g,h\in \Gamma,
\end{equation}
but there is no sequence of homomorphisms $\sigma_k:\Gamma\to U(n_k)$ such that 
\begin{equation}\label{eq:2norm}
\lim_{k\rightarrow\infty} \|  \pi_k(g)-\sigma_k(h) \|_2 =0 \quad \text{ for all }g\in \Gamma.
\end{equation}
In particular this applies to the groups $G_p,K_p$ and $\tilde{G}_p$ of this section.
\end{theorem}
\begin{proof}  By 
\cite[Corollary C]{Tikuisis17} the canonical trace on $\Gamma/\Delta$ is quasidiagonal \footnote{Technically \cite[Corollary C]{Tikuisis17} only states that $C^*_r(\Gamma/\Delta)$ is quasidiagonal, but the proof shows that the canonical trace is quasidiagonal.}.  This means there is a sequence of  functions $\tilde{\pi}_k:\Gamma \to M_{n_k}(\C)$ such that
\begin{enumerate}
\item Equation (\ref{eq:almostmap}) holds.
\item $\tilde{\pi}_k(g) = 1$ for all  $g\in \Delta.$
\item $\tilde{\pi}_k(g^{-1}) = \tilde{\pi}_k(g)^*$ for all $g\in \Gamma$
\item $\|\tilde{\pi}_k(g)\|_\infty \leq 1$ for all $g\in \Gamma$
\item $\lim_k \tr_{n_k}\circ\tilde{\pi}_k(g) = \left\{ \begin{array}{cc} 1 & \text{ if } g\in \Delta\\
                                                                                0 & \text{ if } g\not\in \Delta \end{array} \right.$
 \end{enumerate}
For each $k$ and $g\in \Delta$ we define $\pi_k(g)$ to be the unitary in the polar decomposition of $\tilde{\pi}_k(g)$ if $\tilde{\pi}_k(g)$ is invertible and the identity matrix if $\tilde{\pi}_k(g)$ is not invertible. By conditions (1),(2) (applied to the identity) and (3)  we have
\begin{equation*}
\lim_{k\rightarrow\infty} \| \tilde{\pi}_k(g)\tilde{\pi}_k(g)^* - 1  \|_\infty =0.
\end{equation*} 
Hence by Lemma \ref{lem:basicalc} the sequence of maps $\pi_k$ still satisfy equation (\ref{eq:almostmap})  and condition (5) above but they now map into $U(n_k).$

Suppose there is a sequence $\sigma_k$ satisfying (\ref{eq:2norm}) then by the Cauchy-Schwarz inequality we have
\begin{equation*}
\lim_{k\rightarrow\infty}|\tr_{n_k}\circ \pi_k(g) -  \tr_{n_k}\circ \sigma_k(g)| \leq  \lim_{k\rightarrow\infty}\|\pi_k(g) - \sigma_k(g)\|_2=0,
\end{equation*}
for all $g\in \Gamma.$ Hence by (5) above the trace $\tau$ is a pointwise limit of finite-dimensional traces, a contradiction.
\end{proof}

\begin{remark} 

Let $\tau$ be as in Theorem \ref{thm:sform}.  The C*-algebra generated by the GNS representation of $\tau$ is the group C*-algebra $C^*_r(\Gamma/\Delta)$ which satisfies the UCT and allowed us to use  \cite[Corollary C]{Tikuisis17}.  To apply the method of proof of Theorem \ref{thm:sform} to another ``bad" pair of group $\Gamma$ and trace $\phi$ one would first need to verify that the C*-algebra generated by the GNS representation of $\phi$ satisfies the UCT.  It is unknown if all C*-algebras generated by representations of amenable groups satisfy the UCT.
\end{remark}
\begin{remark} An extremely motivated individual could possibly get around using the deep results of \cite{Tikuisis17} and all of the deep operator algebraic results it relies on. One could perhaps adapt the projection building techniques of Orfanos \cite{Orfanos11} for residually finite groups to build explicit maps $\pi_k$ from Theorem \ref{thm:sform}.  One cannot directly apply his techniques because the quotient $\Gamma/\Delta$ is not residually finite in any of our examples.  
\end{remark}
\bibliographystyle{alpha}

\end{document}